\documentclass[12pt,fleqn]{amsart}

\usepackage{amsrefs}
\usepackage{bookmark}
\usepackage{graphicx} 
\usepackage{enumerate}
\usepackage{korpi-lib}
\usepackage{tikz-cd}
\usepackage{comment}

\title{Barely alternating real almost chains and extension operators for compact lines}
\author[A. Aviles]{Antonio Avilés}
\author[M. Korpalski]{Maciej Korpalski}
\address{Universidad de Murcia, Departamento de Matemáticas, Campus de Espinardo 30100 Murcia, Spain.}
\email{avileslo@um.es}
\address{Instytut Matematyczny, Uniwersytet Wroc\l awski, pl.\ Grunwaldzki 2, 50-384 Wroc\-\l aw\\ Poland}
\email{Maciej.Korpalski@math.uni.wroc.pl}
\date{}
\thanks{Research suupported by PID2021-122126NB-C32 funded by MCIN/AEI/10.13039/501100011033 and by “ERDF A way of making Europe”and by Fundación Séneca - Agencia de Ciencia y Tecnología de la Región de Murcia (21955/PI/22)}

\begin{document}

\begin{abstract}
Assume $\MA$. We show that for every real chain of size $\kappa$ in the quotient Boolean algebra $P(\omega)/fin$ we can find an almost chain of representatives such that every $n\in\omega$ oscillates at most three times along the almost chain. This is used to show that for every countable discrete extension of a separable compact line $K$ of weight $\kappa$ there exists an extension operator $E:C(K)\longrightarrow C(L)$ of norm at most three.
\end{abstract}

\maketitle

\section{Introduction} \label{sec:1 Introduction}

A real almost chain is a family of sets $\mathcal{A} = \{A_x : x\in X\}$ indicated in a set $X\subseteq \mathbb{R}$ increasingly with respect to almost inclusion. That is, $A_x\subseteq^* A_y$ when $x<y$, where $R\subseteq^* S$ means that $R\setminus S$ is finite. A first example can be given by $A_x = \{q\in\mathbb{Q} : q<x\}$ for $x\in\mathbb{R}$. This is in fact a chain: $A_x\subset A_y$ when $x<y$. An easy way to convert this into an almost chain that is not a chain could be to add a finite set to each $A_x$. But this would just be a \emph{finite adjustment} of the previous chain, rather than something essentially new.

\begin{definition}
    $\{B_x : x\in X\}$ is said to be a finite adjustment of $\{A_x: x\in X\}$ if $A_x =^* B_x$ for all $x$.
\end{definition}

Here, $A=^* B$ means that $A\subseteq^* B$ and $B\subseteq^* A$. For a second example, given an irrational number $x\in (0,1)$ with binary expression 
$x= 0. x_1 x_2 x_3\cdots$, $x_i\in\{0,1\}$ consider the sequence $$S_x = \{0.x_1\cdots x_n 0 : x_{n+1}=1\}, $$ and define $A'_x = \{q\in\mathbb{Q}_2 : q<x, q\not\in S_x\}$, where $\mathbb{Q}_2$ are the numbers with finite binary expansion. This time, not only this is not a chain for inclusion, one can see that no uncountable subfamily is a finite adjustment of a chain. Examples of this sort due to Marciszewski were considered in \cite{KP23}*{Construction 5.1} in relation with extension operators from compact lines. Motivated by the same kind of problems, we investigate if almost chains can be readjusted to become \emph{barely alternating} rather than chains.

\begin{definition}
We say that $\{A_x : x \in X\}$ is \textit{barely alternating} if $A_{x_1}\cap A_{x_3}\subseteq A_{x_2}\cup A_{x_4}$ whenever $x_1<x_2<x_3<x_4$.
\end{definition}

Our main result regarding this is the following:

\begin{theorem} \label{mr:2 The aligning forcing}
Under $\MA$, for a set $X$ of cardinality $\kappa$, every real almost chain $\{A_x : x \in X\}$ of subsets of $\omega$ has a barely alternating finite adjustment $\{B_x : x \in X\}$.
\end{theorem}

This resembles the phenomenon shown in \cite{AT11} that, although $(\omega_1,\omega_1)$-gaps exists in ZFC by Hausdorff's classical construction, there exist no $(\kappa,\kappa,\kappa)$-triple gaps under Martin's Axiom $\MA$. The main application of Theorem\ref{mr:2 The aligning forcing} will be the following:

\begin{theorem} \label{apl:2 Ext operator of norm <=3}
Under $\MA$, if $K$ is a separable compact line of weight $\kappa$, then for each countable discrete extension $L$ of $K$ there is an extension operator $E: C(K) \to C(L)$ of norm at most three.
\end{theorem}

Let us explain the terminology. A compact line is a compact space $K$, so that there is a linear order on $K$ in which all sets $\{x: x>a\}$ and $\{x: x<b\}$ are open. For a compact space $K$ we say that a compact space $L$ is a \textit{countable discrete extension} if $K \sub L$ and $L \sm K$ is a discrete countable space. We always identify $L \sm K$ to $\omega$ in such a situation. For a compact space $K$ we write $C(K)$ for a Banach space of real-valued continuous functions on $K$ with the supremum norm. 
When we consider two compact spaces $K \sub L$ we say that a bounded linear operator $E: C(K) \to C(L)$ is an \textit{extension operator} if for all $f \in C(K)$ we have $Ef|_K = f$.\\

Theorem~\ref{apl:2 Ext operator of norm <=3} solves the problem posed at the end of the introduction in \cite{KP23}. 
In that paper it has been shown that if the weight of a separable compact line $K$ is large enough (the cardinality of a subset of reals that cannot be covered by a sequence of closed sets of measure zero), 
then $K$ has a countable discrete extension for which there is no extension operator $E:C(K)\longrightarrow C(L)$. 
Our result proves that just assuming uncountable weight is consistently not enough. If $K$ has a countable weight, an extension operator of norm one exists by the classical result of Dugundji, cf. \cite{PE68}*{Theorem 6.6}. 
When it exists, the minimal norm of an extension operator from a compact line to a countable discrete extension is always an odd number \cite{KP23}. 
It remains an open problem whether an analogous result of Theorem~\ref{apl:2 Ext operator of norm <=3} holds for general separable compact spaces instead of compact lines, cf. \cite{KP23}*{Problem 7.1} and \cites{AMP20,CT19,MP18} for context and the relation with twisted sums of $C(K)$ spaces and $c_0$.\\

For the remaining of the paper, Section~\ref{sec:2 Main forcing} will be devoted to the proof of Theorem~\ref{mr:2 The aligning forcing}, while Section~\ref{sec:3 Applications to CL} contains the proof of Theorem~\ref{apl:2 Ext operator of norm <=3}.

\section{Barely alternating almost chains} \label{sec:2 Main forcing}

For strictly technical reasons we recall here one classical fact from the theory of gaps. For a reference, see e.g. \cite{TJ03}*{Lemma 29.6.}. 

\begin{lemma} \label{mr:1 The omega-omega gaps}
There are no $(\omega, \omega)$ gaps, i.e. if $\{U_n \sub \omega: n\in \omega\}$ is an ascending almost chain and $\{V_n \sub \omega: n\in \omega\}$ is a descending almost chain such that for all $n, m \in \omega$ we have $U_n \sub^* V_m$, then there is a set $U \sub \omega$ such that $U_n \sub^* U \sub^* V_n$ for all $n \in \omega$.
\end{lemma}

We will follow the convention in forcing that a stronger condition is greater, as in \cite{DF84}. Entering into the proof of Theorem~\ref{mr:2 The aligning forcing}, without loss of generality we can assume that $X$ is a dense subset of $[0,1]$. As $\mathbb{R}$ is order-isomorphic to $(0,1)$ we can reindex inside $(0,1)$, and then we could increase it by a countable dense set. The almost chain $\cA$ could be then extended on this countable set using Lemma \ref{mr:1 The omega-omega gaps} (as $X \sub [0, 1]$, all gaps have to be at most $(\omega, \omega)$).
Let us consider the following partial order:
\begin{align*}
\bP = \{(F, \cB = \{B_x \sub \omega: x \in F\}) : \ 
& F \sub X, F \text{ is finite}, \\
& A_x =^* B_x \text{ for } x \in F, \\
&\cB \text{ is barely alternating}
\}, \\
(F_1, \cB_1) \leq (F_2, \cB_2) &\iff F_1 \sub F_2 \land \cB_1 \sub \cB_2.
\end{align*}
We will prove that:
\begin{enumerate}[a)]
    \item $\bP$ is ccc;
    \item by Martin's Axiom, we can find an almost chain as in the statement of the theorem.
\end{enumerate}

a) Consider any uncountable set $P = \{(F^\alpha, \cB^\alpha)\}$ of conditions from $\bP$. Without loss of generality we can assume that (as $P$ is uncountable and we can pass to uncountable subfamilies):
\begin{enumerate}[$\bullet$]
    \item Sets $F^\alpha$ form a $\Delta$-system, so for all $\alpha \in \kappa$ they are of the form $F^\alpha = R \cup G_\alpha$, with $G_\alpha$ being pairwise disjoint sets of fixed cardinality $k\in \omega$. 
    For each $\alpha \in \kappa$ let $G_\alpha = \{g^\alpha_i : 1 \le i \le k\}, R = \{r_i : 1 \le i \le l\}$ be increasing enumerations of those sets. 
    We can also fix for each $\alpha \in \kappa$ and $i \in \{1, \dots, l-1\}$ the amount of points in $G_\alpha$ below $r_1$, between $r_i$ and $r_{i+1}$, and above $r_l$.
    \item For each $1 \leq i \leq k-1$ there are some rational points $p_i, q_i$ satisfying $g^\alpha_i < p_i < q_i < g^\alpha_{i+1}$ for every $\alpha \in \kappa$.
    \item There is a set $Z \sub X$ such that $R \sub Z$ and between each points $g^\alpha_i$ and $g^\alpha_{i+1}$ there is a point $z_i \in Z$ chosen either from $R$ (if there is any; there can be multiple of them) or from the interval $(p_i, q_i)$. Write $Z = \{z_j : 1 \le j \le k-1\}, z_i < z_j$ for $i < j$. 
    \item As the set $R$ is finite, it is possible to fix the set $B_r$ for all $r\in R$ and $\alpha < \kappa$.
    We can also choose sets $B_z =^* A_z$ for $z \in Z\sm R$ such that the almost chain $\cB^\alpha \cup \{B_z : z\in Z\sm R\}$ is still barely alternating.
    \item There is a finite set $D \sub \omega$ such that, for all $\alpha \in \kappa$, the almost chain $\cC = \{B^\alpha_{g^\alpha_i} : 1 \le i \le k\} \cup \{C_z : z \in Z\}$ is a chain outside of the set $D$, i.e. for all $A, B \in \cC$ if $A \sub^*B$, then $A \sm B \sub D$.
    For each $d\in D$ we can also fix for which $i$~the set $B^\alpha_{g^\alpha_i}$ contains the element $d$, so that for all $\alpha, \beta \in \kappa$ we have $d \in B^\alpha_{g^\alpha_i} \iff d \in B^\beta_{g^\beta_i}$
\end{enumerate}

Now let us take two elements $(F^\alpha, \cB^\alpha), (F^\beta, \cB^\beta) \in \bP$ and check whether they are compatible. It is enough to show that the almost chain $\cB^\alpha \cup \cB^\beta$ is barely alternating. Write as $B_x$ the element of $\cB^\alpha \cup \cB^\beta$ with the index $x\in F^\alpha \cup F^\beta$.

Fix any $m \in \omega$ and elements $x_1 < x_2 < x_3 < x_4$ from $F^\alpha \cup F^\beta$. If $m \in D$, then as $m \in B^\alpha_{g^\alpha_i} \iff m \in B^\beta_{g^\beta_i}$ for all $i$ and both almost chains $\cB^\alpha, \cB^\beta$ are barely alternating, we have $m\in B_{x_1}\cap B_{x_3} \implies m \in B_{x_2}\cup B_{x_4}$.

For $m \notin D$, if $m \in B_{x_1}$ and $m\notin B_{x_2}$, by the definition of set $D$ we have $x_1, x_2 \notin R \sub Z$. Without loss of generality we can thus denote $x_1 = g^\alpha_i, x_2 = g^\beta_j$ for some $1 \leq i \leq j \leq k$.
Now if $i < j$, then we have again a contradiction with the definition of $D$, as $m$ has to be an element of $B_{z_i}$ (as it is an element of $B^\alpha_{g^\alpha_i}$) and cannot be an element of $B_{z_i}$ (as it is not an element of $B^\beta_{g^\beta_j}$).
If $i = j$, then $m\in B_{z_i}$ and $m\in B_{x_3}, m\in B_{x_4}$, as $x_3, x_4 \geq z_i \in Z$, so the almost chain $\cB^\alpha \cup \cB^\beta$ is indeed barely alternating. \\

b) We want to prove that for all $x \in X$ the set $\cD_x = \{(F, \cB) \in \bP : x \in F\}$ is dense in $\bP$. As $|X| = \kappa$, from $\MA$ we will have a filter $\cG \sub \bP$ which intersects all sets $\cD_x$, as $\bP$ is a ccc partial order. Then the family $\bigcup_{(F, \cB) \in \cG} \cB$ is the almost chain we are looking for.

Fix any $x\in X$ and an element $(F, \cB) \in \bP$. If $x\in F$, we are done. If not, let $y, z \in F$ be elements directly preceding and succeeding $x$ in $F$ and put $A = B_y, C = B_z$ (if $x$ is below or above all elements of $F$, put $A = \emptyset$ or $C = \omega$ accordingly). 
Let us put $B_x = (A_x \cup A) \sm (A_x \sm C)$. We can see that $B_x =^* A_x$ as we have added or removed only a finite amount of elements. 

Now we want to check that $(F \cup \{x\}, \cB \cup \{B_x\})$ is an element of $\bP$, so we have to verify that it barely alternating. For all $m\in \omega$ we have either $(m \in B_x \iff m \in A)$ or $(m \in B_x \iff m \in C)$, so the chain $\cB \cup \{B_x\}$ has to be barely alternating. 
It means for all $x\in X$, the set $\cD_x$ is dense in $\bP$ and the proof is complete.

\section{Countable discrete extensions of compact lines} \label{sec:3 Applications to CL}

We embark now in the proof of Theorem~\ref{apl:2 Ext operator of norm <=3}. We will begin by stating some elementary facts about compact lines that will be useful to us. Compact lines can be viewed as linear orders in which all subsets have a supremum (equivalently an infimum) endowed with the order topology. One can find in \cite{EN89}*{3.12.3-4} and \cite{SS78}*{II.39}. Given a compact line $K$, a point $x\in K$ is isolated from the left if $x\not\in\overline{\{y\in K : y<x\}}$, and isolated from the right if $x\not\in\overline{\{y\in K : y>x\}}$. Thus, $K$ can be decomposed into the following sets:
\begin{enumerate}
    \item $K_{\leftarrow}$ is the set of points that are isolated from the left but not from the right.
    \item $K_{\rightarrow}$ is the set of points that are isolated from the right but not from the left.
    \item $K_{\bullet}$ is the set of isolated points.
    \item $K_{-}$ is the set of points that are not isolated from any side.    
\end{enumerate}

\begin{proposition}\label{cardinalinvariants}
    $weight(K)  = \max\{dens(K), |K_{\rightarrow}|\}$
\end{proposition}
\begin{proof}
   For the inequality $[\geq]$, we know that $dens(K)\leq weight(K)$, and if $\mathcal{B}$ is a basis of open sets of $K$, then for any $x\in K_{\rightarrow}$, there should be $W_x\in \mathcal{B}$ such that $x\in W_x\subseteq \{y : y\leq x\}$. Since $\max(W_x)=x$, this shows that $|\mathcal{B}|\geq |K_{\rightarrow}|$. For the inequality $[\leq]$, fix a dense subset $D$ of $K$. Then, the family of intervals
   $$\left\{(a,b),[\min(K),b),(a,x], \{p\} : a,b\in D\cup K_{\leftarrow}\cup K_{\bullet}, x\in K_\rightarrow, p\in K_{\bullet} \right\}$$
   forms a basis for the topology. Since $|K_\bullet|\leq dens(K)$, we are done.
\end{proof}

\begin{proposition}\label{density}
    The set $K\setminus K_\leftarrow$ is dense in $K$.
\end{proposition}

\begin{proof}
    If not we would have a nonempty open interval $(a,b)\subset K_\leftarrow$. In particular, this interval does not contain isolated points and whenever $x,y\in (a,b)$ with $x<y$ there exists $z\in (x,y)$. This allows to construct points $x_1<x_2<\cdots<y$ in $(a,b)$ and then $\sup\{x_1,x_2,\ldots\} \in (a,b)\setminus K_\leftarrow$, a contradiction.
\end{proof}

The following is a reformulation of the classical Ostaszewski's result \cite{OS74} that every separable compact line can be obtained by \emph{splitting} the points of compact subset of $[0,1]$. 

\begin{proposition}\label{realorder}
    If $K$ is a separable compact line, then $K\setminus K_\leftarrow$ is order-isomorphic to a subset of the real line.
\end{proposition}
\begin{proof}
Let $D$ be a countable dense subset of $K$. Since the lexicographic product $D\times_{lex} \{0,1\}$ is a countable linear order, there is a strictly increasing function $u:D\times_{lex}\{0,1\} \longrightarrow \mathbb{Q}\cap [0,1]$. Send each $y\in K\setminus K_\leftarrow$ to any real number $r(y)$ satisfying
$u(d,1) \leq r(y) \leq u(d',0)$ whenever $d\leq y < d'$. This is strictly increasing: If $y_1<y_2$ and the interval $(y_1,y_2)$ is nonempty, picking $d\in (y_1,y_2)\cap D$ we get $r(y_1)\leq u(d,0)<u(d,1)\leq r(y_2)$. If that interval is empty, then $y_2$ is isolated from the left, hence also from the right since $y_2\not\in K_{\leftarrow}$. That implies that $y_2\in D$ and then $r(y_1)\leq u(y_2,0)<u(y_2,1) \leq r(y_2)$.
\end{proof}

By Propositions \ref{cardinalinvariants} and \ref{density}, we can fix a dense subset $Y$ of $K$  of cardinality the weight of $K$ and such that 
$$K_{\rightarrow} \cup K_{\bullet} \subseteq Y \subseteq K \setminus K_{\leftarrow}.$$ 
Moreover, when $K$ is separable, $Y$ will be order-isomorphic to a subset of the real line by Proposition \ref{realorder}.

We fix a countable discrete extension $L$ of the separable compact line $K$.

\begin{lemma}
  For every $y\in Y$ there is a continuous function $h_y\in C(L)$ such that \begin{itemize} \item $h_y(x) < 0$ if $x<y$, \item $h_y(x) > 0$ if $x>y$, \item $h_y(y)<0$ if $y$ is isolated from the right.\end{itemize}
\end{lemma}

\begin{proof}
    By Tietze's extension theorem, it is enough to define $h_y$ on $K$. If $y$ is isolated from the right we can define $h_y(x)=-1$ if $x\leq y$ and $h_y(x)=1$ if $x> y$. Otherwise, by the definition of $Y$, $y$ is not isolated neither from left nor right. By separability, this means that there exists an increasing sequence $(a_n)$ and a decreasing sequence $(b_n)$ in $K$ that both converge to $y$. This implies that $\{x\in K : x<y\}$ and $\{x\in K : x>y\}$ are disjoint $F_\sigma$ open subsets of $K$ and hence we can define a continuous function that is negative on one set and positive on the other. 
\end{proof}

We fix functions $h_y$ as in the previous lemma. We define
$$A_y = \{n<\omega : h_y(n)<0\}.$$
This is an almost chain of subsets of $\omega$. Indeed, if $y<y'$, then every $n\in A_y\setminus A_{y'}$ satisfies $h_y(n)<0$ and $h_{y'}(n)\geq 0$. Therefore
$$K \cap \overline{A_y\setminus A_{y'}} \subseteq \{x\in K : h_y(x)\leq 0, h_{y'}(x)\geq 0\}\subseteq\{x\in K : y'\leq x\leq y\} = \emptyset,$$
hence $A_y\setminus A_{y'}$ is finite. 

Using Theorem \ref{mr:2 The aligning forcing}, we can finitely adjust the almost chain $\{A_y : y \in Y\}$ into a barely alternating almost chain $\{B_y : y \in Y\}$. Then for each $n \in \omega$ define in $K$ points $x^0_n, x^1_n, x^2_n$ as
\[x^0_n = \inf \{y \in Y : n \in B_y\},\]
\[x^1_n = \inf \{y \in Y : n \notin B_y \text{ and } y > x^0_n\},\]
\[x^2_n = \inf \{x \in Y : n \in B_y \text{ and } y > x^1_n\}.\]
Notice that, formally, $\inf(\emptyset) = \max(K)$. Using the points $x^0_n, x^1_n, x^2_n$ we can define an extension operator $E: C(K) \to C(L)$. For $f\in C(K)$ write $Ef|_K = f$ and for $n\in \omega = L\sm K$ put $Ef(n) = f(x_n^0) - f(x_n^1) + f(x_n^2).$

It is clear that $E$ is linear and that its norm will be at most 3. We only need to verify that $Ef$ is continuous on $L$ for all $f\in C(K)$. Since $L$ is first-countable the topology is determined by sequences, so it is enough to check that if a sequence $(n_k)_{k\in \omega}\subset L\setminus K$ converges to $z\in K$, then $\lim_k Ef(n_k) = f(z)$. If this was not true, there would be a subsequence with $\lim_i Ef(n_{k_i}) = \ell \neq  f(z)$. Passing to a subsequence, we can assume that each of the sequences $(x^0_{n_k})$, $(x^1_{n_k})$ and $(x^2_{n_k})$ is monotone converging to points $x^0\leq x^1 \leq x^2$ in $K$.\\

If $y\in Y$, then we have the following:

\begin{enumerate}
    \item If $y< x^0$, then for all but finitely many $k$, we will have $y<x^0_{n_k}$, hence $n_k\not\in B_y$. Therefore for all but finitely many $k$ we will have $n_k\not\in A_y$, hence $h_y(n_k)\geq 0$. Therefore $h_y(z)\geq 0$, so $y\leq z$.

    \item If $x^0<y<x^1$, then for all but finitely many $k$, we will have $x^0_{n_k}
    <y<x^1_{n_k}$, hence $n_k\in B_y$. Therefore for all but finitely many $k$ we will have $n_k\in A_y$, hence $h_y(n_k)\leq 0$. Therefore $h_y(z)\leq 0$, so $y\geq z$.

    \item If $x^1<y<x^2$, then for all but finitely many $k$, we will have $x^1_{n_k}
    <y<x^2_{n_k}$, hence $n_k\not\in B_y$. Therefore for all but finitely many $k$ we will have $n_k\not\in A_y$, hence $h_y(n_k)\geq 0$. Therefore $h_y(z)\geq 0$, so $y\leq z$.

     \item If $x^2<y$, then for all but finitely many $k$, we will have $x^2_{n_k}
    <y$, hence $n_k\in B_y$. Therefore for all but finitely many $k$ we will have $n_k\in A_y$, hence $h_y(n_k)\leq 0$. Therefore $h_y(z)\leq 0$, so $y\geq z$.
\end{enumerate}

Since $Y$ is dense in $K$,  conditions (1) and (4) imply that $x^0 \leq z \leq x^2$.\\

Claim 0. If $x^0 < x^1$, then $z= x^0$. To see this, we distinguish two cases. If the interval $(x^0,x^1)$ is nonempty, then we just use (2) and the density of $Y$. If $(x^0,x^1)=\emptyset$ then we have $x^0\in K_\rightarrow \cup K_{\bullet}\subset Y$. Also $x^0_{n_k}\leq x^0 < x^1 \leq x^1_{n_k}$ for all but finitely many $k$. If the former inequalities are strict, we can reason as in (2) and get that $h_{x^0}(n_k)\leq 0$ for infinitely many $k$ and conclude that $h_{x^0}(z)\leq 0$ so $z\leq x^0$, as desired. The other case is that $x^0 = x^0_{n_k}$ for infinitely many $k$. Since $x^0$ has an immediate successor, it has to be a minimum when it is an infimum,
        
        $x^0 = \min\{y\in Y: n_k\in B_y\}$ for all $k$,

        so $n_k\in A_{x^0}$ (so $h_{x^0}(n_k)\leq 0$) for infinitely many $k$, hence $h_{x^0}(z)\leq 0$ and $x^0\geq z$ again.\\

Claim 1. If $x^1<x^2$, then $z=x^2$. If the interval $(x^1,x^2)$ is nonempty, then we just use (3) and the density of $Y$. If $(x^1,x^2)=\emptyset$ then we have $x^1\in K_\rightarrow \cup K_{\bullet}\subset Y$. Also $x^1_{n_k}\leq x^1 < x^2 \leq x^2_{n_k}$ for all but finitely many $k$. If the former inequalities are strict, we can reason as in (3) and get that $h_{x^1}(n_k)\geq 0$ for infinitely many $k$ and conclude that $h_{x^1}(z)\geq 0$. Since $x^1$ is isolated from the right, we conclude that $x^1 < z$, therefore $x^2\leq z$ as desired. The other case is that $x^1 = x^1_{n_k}$ for infinitely many $k$. Since $x^1$ has an immediate successor, it has be to a minimum when it is an infimum,
        
        $x^1 = \min\{y\in Y: n_k\not\in B_y, y>x^0_{n_k}\}$ for all $k$,

        so $n_k\not\in A_{x^1}$ (so $h_{x^1}(n_k)\geq 0$) for infinitely many $k$, hence $h_{x^1}(z)\geq 0$ and $x^2\leq z$ again.\\ 

From the two claims, the situation $x^0<x^1<x^2$ is impossible and we must have one of the following three:
\begin{itemize}
    \item $x^0 < x^1 = x^2$ and $z=x^0$,
    \item $x^0 = x^1 < x^2$ and $z=x^2$,
    \item $x^0 = x^1 = x^2 = z$.
\end{itemize}
In any case,$$\lim_k Ef(n_k) = \lim_k f(x^0_{n_k}) - f(x^1_{n_k}) + f(x^2_{n_k}) = f(x^0) - f(x^1) + f(x^2) = f(z)$$ as desired.

\bibliography{refs}

\begin{bibdiv}
\begin{biblist}

\bib{AMP20}{article}{
      author={Avilés, A.},
      author={Marciszewski, W.},
      author={Plebanek, G.},
       title={Twisted sums of $c_0$ and {C(K)}-spaces: A solution to the {CCKY} problem},
        date={2020},
     journal={Advances in Mathematics},
      volume={369},
       pages={107168},
}

\bib{AT11}{article}{
      author={Avilés, A.},
      author={Todorcevic, S.},
       title={Multiple gaps},
        date={2011},
     journal={Fund. Math.},
      volume={213},
       pages={15\ndash 42},
}

\bib{CT19}{article}{
      author={Correa, C.},
      author={Tausk, D.~V.},
       title={Local extension property for finite height spaces},
        date={2019},
     journal={Fundam. Math.},
      volume={245},
      number={2},
       pages={149\ndash 165},
}

\bib{EN89}{book}{
      author={Engelking, R.},
       title={General {T}opology},
   publisher={Heldermann Verlag, Berlin},
        date={1989},
}

\bib{DF84}{book}{
      author={Fremlin, D.H.},
       title={Consequences of {M}artin's {A}xiom},
      series={Cambridge Tracts in Mathematics},
   publisher={Cambridge University Press},
        date={1984},
}

\bib{TJ03}{book}{
      author={Jech, T.},
       title={Set {T}heory: {T}he {T}hird {M}illennium {E}dition},
   publisher={Springer},
        date={2003},
}

\bib{KP23}{unpublished}{
      author={Korpalski, M.},
      author={Plebanek, G.},
       title={Countable discrete extensions of compact lines},
        date={2023},
        note={arXiv:2305.04565},
}

\bib{MP18}{article}{
      author={Marciszewski, W.},
      author={Plebanek, G.},
       title={Extension operators and twisted sums of $c_0$ and {C(K)} spaces},
        date={2018},
     journal={Journal of Functional Analysis},
      volume={274},
       pages={1491\ndash 1529},
}

\bib{OS74}{article}{
      author={Ostaszewski, A.J.},
       title={A characterization of compact, separable, ordered spaces},
        date={1974},
     journal={J. London Math. Soc.},
      volume={7},
       pages={758–760},
}

\bib{PE68}{article}{
      author={Pe{\l}czy{\'n}ski, A.},
       title={Linear extensions, linear averagings, and their applications to linear topological classification of spaces of continuous functions},
        date={1968},
     journal={Diss. Math.},
      volume={58},
}

\bib{SS78}{book}{
      author={Steen, L.~A.},
      author={Seebach, J.~A.},
       title={Counterexamples in topology. 2nd ed},
   publisher={{Springer}-{Verlag}},
        date={1978},
}

\end{biblist}
\end{bibdiv}

\end{document}